\documentclass[12pt]{amsart}
\usepackage{amssymb}
\theoremstyle{plain}
\newtheorem{theorem}{Theorem}
\newtheorem{proposition}{Proposition}

\theoremstyle{definition}

\theoremstyle{remark}

\numberwithin{equation}{section}
\begin{document}
\title[Analyticity for non-linear sums of squares]{Local Real Analyticity
of Solutions for\\ sums
of squares of non-linear vector fields}
\author{David S. Tartakoff}
\address{Department of Mathematics, Statistics and Computer Science\\
					University of Illinois at Chicago \\
					851 So. Morgan St., Chicago
IL  60607 USA}
\email{dst@uic.edu}
\author{Luisa Zanghirati}
\address{Dipartimento di Matematica, Universit\`a di Ferrara, via
Machiavelli 35,
44100 Ferrara, Italia.}
\email{zan@dns.unife.it}

\begin{abstract}
We show that all smooth solutions of model non-linear sums of squares of vector
fields are locally real analytic.
\end{abstract}
\maketitle

\section{Introduction}

   We consider sums of squares of non linear vector fields,
that is equations
such as
$$P(x,u,D)u=\sum_1^qX_j^2u=f$$
with the new feature that $\{X_j\}$ may depend in their
``coefficients'' on the solution $u.$
As a prime example of this class we consider the following case (for
$r>0$):
\begin{equation} \label{1}
P_u(D)v:=\left((D_x)^2+(x^rD_t)^2+(x^r\tilde h(x,t,u)D_t)^2\right)v,\
(x,t)\in\mathbb{R}^2\end{equation}  with $\tilde h$ real valued and real
analytic in its arguments.

We shall assume our solution $u$ to be $C^\infty,$ 
since smoothness (starting from $C^{2+\alpha}$)  follows from the
  arguments of Xu \cite{Xu1} which are based on the
subelliptic estimate clearly satisfied by $P_u$ and the paradifferential
calculus of Bony \cite{B}.

\section{Results}
\begin{theorem} If $f$ is real analytic near $(x_0,t_0),$
then so is any smooth solution to (\ref{1}).
\end{theorem}

We remark that the problem is significant in its own right
and also because
it bears the same resemblance to general quasilinear subelliptic
partial differential
equations that the sums of squares of linear vector fields do to the
subelliptic complexes and `boundary Laplacians'
arising from the $\overline{\partial}_b$ operator
in several complex variables. In particular, the local real analytic
hypoellipticity of those (in the linear case) with symplectic
characteristic variety (roughly corresponding to $r=1$ here),  proved
independently by Treves and Tartakoff in 1978 (\cite{Ta78},
\cite{Ta81}, \cite{Tr78}), propels one quite reasonably to ask the
same question in the quasilinear setting, of which the type of
operator under study here is a simple prototype. (NB - the vector
fields arising from
$\overline{\partial}_{(b)}$ correspond more directly to
$\partial_t-y\partial_t$ and $\partial_y+x\partial_t$ than to
$\partial_t,\partial_y, x\partial_t$ and $y\partial_t$ as separate vector
fields; nonetheless the ``Grushin-type'' operators have always provided
the most tractable models.)

\section{Proof}

   Using standard arguments it is easy to prove the following
{\em a priori} estimates: $\forall s\geq 0, u\in C^\infty$ and compact
${\mathcal U}, \exists C=C_{s,u, {\mathcal U}}: \forall\ v\in
C_0^{\infty}({\mathcal U}),$
\begin{equation} \label{se1}
{\sum_1^3\|X_iv\|_s^2+\|v\|^2_{s+\frac1{r+1}}\leq
C\bigl\{|(P_uv,v)_s|+\|v\|_s^2\bigr\}\ } {\text and}
\end{equation}
\begin{equation}\label{se2}
\sum_{i,j=1}^3\|X_iX_jv\|_s^2+
\sum_1^3\|X_iv\|^2_{s+\frac1{r+1}}+\|v\|^2_{s+\frac2{r+1}}\leq
C\bigl\{\|P_uv\|_s^2+\|v\|_s^2\bigr\}
\end{equation}
 where
$\|\cdot\|_s=\|\cdot\|_{H^s},\ P_u\equiv P(x,u,D),\,
X_1=D_t,\ X_2=x^rD_t,\ X_3=X_3^{(u)}=x^r\tilde h(t,x,u)D_t,$ and
$C$ depends only on the first $s+3$ derivatives of $u.$ 

However the estimate we will need uses the maximality (and arbitrary
positivity) of (\ref{se2}) rather than its subellipticity: 
with $|||v|||_s$ defined as follows for $s$ a positive integer, 
\begin{equation}
|||v|||_s \equiv \sum_{|I| \leq 2}\|X^Iv\|_s, \quad 
(|||v|||_{H^s({\mathcal U})} \equiv \sum_{|I| \leq
2}\|X^Iv\|_{H^s({\mathcal U})}, {\text{for }}{\mathcal U} {\text{ open}})
\end{equation}
then for any $K$ and with $X^I = X^{I_1}X^{I_2}\ldots
X^{I_{|I|}}, \,\exists C_K:\forall v \in C_0^\infty({\mathcal U}),$  
\begin{equation}
\label{max}
|||v|||_2 + K \sum_{|I| \leq 1}\|X^Iv\|_2 \leq C\|P_uv\|_2 +
C_K\|v\|_0.
\end{equation}
\section{The general scheme}
The general scheme, as
always, will be to use the {\em a priori} estimate applied to functions
$v=\varphi D^mu$ and then to bring
$\varphi D^m$ to the left of
$P_u$ modulo errors which are handled inductively. Noting that the {\em a
priori} estimate provides for
maximal control (i.e. no loss of derivatives) in the $D_x$ direction,
we limit ourselves to estimating $\varphi D_t^mu$.
$\varphi$ will be a smooth localizing function, namely identically
equal to one in a fixed open set ${\mathcal U}_0,$
where we wish to prove the solution $u$ is analytic, and supported in
${\mathcal U}_1,$ the
open set where the data are assumed to be real analytic. The localizing
function $\varphi(x,t)$ may be taken to be of the form
$\tilde\varphi(t)\tilde{\varphi}(x),$ and terms with derivatives on
$\tilde\varphi(x)$ may be disregarded since the operator is elliptic
when $x$ is away from $0,$ namely in the support of derivatives of
$\tilde\varphi(x).$ . Thus for our purposes, $\varphi =
\varphi(t)$ alone.

Taking  $Pu=0$ without loss of generality, we have from (\ref{max}):

$$
{|||\varphi D_t^mu|||_2\over m!} \sim {\|X_j^2\varphi D_t^mu\|_2\over
m!}+ ...
\leq
{\|P\varphi D_t^mu\|_2\over m!} + ... \leq {\sum\|[X_j^2,\varphi
D_t^m]u\|_2\over m!}+...
$$
$$    \lesssim
\sum_{k=1}^2{\|g_k(x,t,u,u')\varphi^{(k)}x^{2r}D_t^{m+2-k}u\|_2\over
m!}
      +C{\|\varphi x^{2r}[h(u),D_t^mu]D^2_tu\|_2\over m!}+...
$$
\begin{equation}
\lesssim C\sum_{k=1}^2{\|\varphi^{(k)}X^2D_t^{m-k}u\|_2\over m!}+
C{\|\varphi
x^{2r}[h(u),D_t^mu]D_t^2u\|_2\over m!}+...\label{3} \end{equation}
considering $x^{2r}D_t^2=X^2,$ writing $h(\cdot)\equiv\tilde
h^2(\cdot)$ and estimating the norm $\|g_k(x,t,u,u')\|_{H^2({\mathcal
U}_1)}$ by a constant. Here the $g_k(x,t,u,u')$ stand for the coefficients,
aside from
$x^{2r},$ which enter when $\varphi$ is differentiated once or twice, and
the dots ``..." denote terms arising from lower order terms in the
operator
$P,$ terms containing fewer $X$'s.

  We focus on the bracket in the 
last norm, the crucial one. To expand $(D_t^{m'}h)(u(x,t)),$ we will
need to use the Fa\`a di Bruno formula or rather, what will suffice, and
probably be more transparent, crude bounds for the results: writing
$$D^m_t g(u(x,t)) = (u'D_u+D_t)^{m-1}u'D_u g(u(x,t));$$
with primes on $u$ denoting $t$ derivatives, writing this roughly as 
$$D^m_t g(u(x,t)) = ((u'\sigma +D_t)^{m-1}u'_{|\sigma = D_g})g,$$
i.e., $\sigma$ becomes a `countem' for the number of derivatives
received by $g.$ Then this is {\it at worst} 
\begin{equation}\label{approxFdB}\sum_{m'}{m \choose
m'}g^{(m-m')}(D_t^{m'}{u'}^{m-m'}).\end{equation}  
Finally, distributing $\alpha$ objects into $\beta$ positions  yields
$$ {D^a {u'}^b \over a!} = \sum_{a_1+\cdots + a_b=a} 
{{u'}^{(a_1)}\over a_1!}\cdots {{u'}^{(a_b)}\over a_b!}$$

Thus we have 
\begin{equation}\label{bracketwithD^m}
{[D_t^m,P]\over m!}=\sum_{m'=1}^{m}
{(D_t^{m'}h(u(x,t)))\over m'!}
\;{(x^{2r}D_t^2)D_t^{m-m'}\over (m-m')!}.
\end{equation}
with (cf. (\ref{approxFdB})):
$${(D_t^{m'}) h(u(x,t))\over m'!}
\sim
\sum_{m'-m''\geq 1}^m{h^{(m'-m'')}\over
{(m'-m'')}!}\;{(D_t^{m''}({u'}^{m'-m''}))\over m''!}$$ 
$$=\sum_{m'-m''\geq 1}^m{h^{(m'-m'')}\over
{(m'-m'')}!}\sum_{\sum_1^{m'-m''} m''_j = m''} 
{D_t^{m''_1}{u'}\over m''_1!}\cdots {D_t^{m''_{m'-m''}}{u'}\over
m''_{m'-m''}!}$$ or in all, with (\ref{3}),  
\begin{equation}
\label{[Tr,P]}
{|||\varphi D_t^mu|||_2\over m!} \lesssim {\|[\varphi D_t^m, P]u\|\over
m!}+\ldots \lesssim 
\sum_{k=1}^2{\|\varphi^{(k)}X^2D_t^{m-k}u\|_2\over m!}+\ldots
\end{equation}
$$+ \sum\left\|\varphi \;{h^{(m'-m'')}\over
{(m'-m'')}!}{D_t^{m''_1}{u'}\over m''_1!}\cdots
{D_t^{m''_{m'-m''}}{u'}\over m''_{m'-m''}!}\;{(
x^{2r}D_t^2)D_t^{m-m'}u\over (m-m')!}\right\|_2
$$
where the sum is over $m\geq m'-m''\geq 1, \sum_{j=1}^{m'-m''}
m''_j = m''$ and so $\sum_{j=1}^{m'-m''} (m''_j+1) =
m'.$

By associating $x^{2r}$ with a different term in the product if necessary,
we may assume that the last term is of greatest order, and hence that the
others are of order {\it at most} $m/2.$

\section{Remarks on the last sum}

Several remarks are in order concerning the last right hand side. 

First of all,
in utilizing the property that $H^2$ is an algebra to take the product of
norms, there will occur a constant raised to the $m'-m''.$ But this is
allowable, since there are $m'-m''$ derivatives on the (analytic) function
$h$ where we expect a constant to that power. 

Secondly, that power that power always corresponds to the increase in
number of terms of the form ${D_t^{m''_j}{u'}/m''_j!}$ inside the norm; in
the end the number of these terms cannot exceed $m,$ hence the constant
cannot exceed $C^m.$ 

Thirdly, we will associate the localizing function $\varphi,$ with the
highest order term and take it out of the norm, introducing another one
which is closely related to the number of derivatives in that term - in this
case $m-m'.$ In bringing $\varphi$ out of the norm there may be one or
two derivatives (or three or four, given the first terms on the right of
(\ref{[Tr,P]})), and while they will presumably balance quite well with
$m!$ we need to be sure that they balance as well with $(m-m')!$ when
$m-m'$ may actually be rather small (a large drop may have occurred all at
once). To this end we make the following observation: as $m$ drops from
$m$ to $m-m',$ there have appeared $m'-m''$ {\it new} lower order terms,
or
$m'-m''+1$ terms of no greater order, counting the principal one. Thus we
have 
\begin{equation}\label{no.ofterms}
(m'-m''+1)(m-m')\geq m;\qquad i.e., {m \over {m-m'}}\leq m'-m''+1,
\end{equation}
the same factor that occurred before, and appears in the number of
derivatives on $h.$ Thus, again, we can afford $(m/m-m')^4$ without
danger.

The fourth observation concerns the effect of the sum. The sum 
corresponds {\it at most} to the number of ways to partition $m$ derivatives
among at most
$m$ functions, generally many fewer. Denoting by $D$ a
derivative ($m$ of them) and by $u$ a copy of $u$ ($t$ of
them) we are faced with the number of ways to `identify' or select
$t$ items (the $u's$) from among $m+t$ items (the $D$'s and
$u$'s) with the understanding that in an expression such as 
\begin{equation}\label{partitions}
\underbrace{\underbrace{DDDDD}_{m_1}u
\underbrace{DDDDD}_{m_2}u
\underbrace{DDDDD}_{m_3}u
\underbrace{DDDDD}_{m_4}u
\cdots
\underbrace{DDDDD}_{m_t} u}_{m \; D's \;{\rm and }\; t(\leq
m)\, 
\;u's}
\end{equation} 
the $D$'s differentiate only the first $u$ following. The answer
is that there are certainly not more than ${m+t \choose
t}\leq 2^{m+t}\leq 2^{2m}=4^m$ ways. And while we have written this out
only for the first complete iteration of the {\it a priori} estimate, it is a
remarkable fact that the form of the sum does not change after multiple
passes, and hence the number of terms involved is subject to the same
bounds. What is more, the same analysis applies after iteration of
(\ref{[Tr,P]}) (cf. below) and thus the sum will also not pose a difficulty in
proving analyticity and may be replaced by a supremum below. 

Finally, when these considerations enter and readability is an issue due to
the length of lines, we shall tacitly replace the sum by a supremum and omit
a constant such as $C^{m'-m''+2}.$ 

\section{The localizing functions and introducing new ones}

\begin{proposition}
For any two open sets $\Omega_0\Subset\Omega_1,$ with
separation $d={\hbox{dist.}}(\Omega_0,\Omega_1^{c})$ and any
natural number $N,$ there exists a universal constant
$C$ depending only on the dimension and a function
$\Psi=\Psi_{\Omega_0,\Omega_1,N}\in C_0^\infty (\Omega_1),
\Psi\equiv 1 \hbox{ on }\Omega_0$ 
with 
\begin{equation}\label{fns:Ehrenpreis}
|D^\beta \Psi| \leq \left({C\over
d}\right)^{|\beta|+1}N^{|\beta|},
\qquad |\beta| \leq 2N,
\end{equation}
\end{proposition}

The first localizing function, $\varphi = \Psi_m,$
satisfies:
\begin{equation}
\label{def:Psi_m} 
\Psi_m\equiv 1  \mbox{ on }\ {\mathcal U}_0\ ,
\Psi_m\in
C_0^{\infty}({\mathcal U}_{1/m})\ ,
\ |\Psi_m^{(k)}|\leq c^k m^{k}, \;k\leq 4,\end{equation}
where we have set, for $a\geq 0:$
\begin{equation}
\label{def:Ua} 
{\mathcal U}_a=\{(x,t)\in{\mathcal U}_1\ :
\mbox{dist}((x,t),{\mathcal U}_0)<a(\text{dist}(\mathcal U_0,
\mathcal U_1^c)\}\ .\end{equation} 

When the first  
localizing function needs to be replaced but, say, 
 $\tilde{m}$ derivatives of $u$ remain to be estimated, 
we shall localize it with a function identically
equal to one on $\mathcal U_{1/m},$ the support of $\Psi_m$
but dropping to zero in a band of width $1/\tilde{m}$ of the
remaining distance ($a(1-1/m)$) to the complement of  $\mathcal
U_1,$ i.e., supported in 
\begin{equation}
\label{supports1}
\mathcal U_{{1\over m}+({1\over
\tilde{m}})(1-{1\over m})}=\mathcal U_{{1\over m}+{m-1\over
m\tilde{m}}}=\mathcal U_{1-(1-{1\over m})(1-{1\over
\tilde{m}})}.
\end{equation}
We shall denote such a function by $_{1\over m}\Psi
_{\tilde{m}}$ 
That is,
${}_\rho\Psi_\sigma$ satisfies:
\begin{equation}
\label{Psirhosigma1} 
{}_\rho\Psi_\sigma \equiv 1 \mbox{ on }
{\mathcal U}_{\rho}, \qquad
{}_\rho\Psi_\sigma \in
C_0^{\infty}({\mathcal{U}}_{\rho + {1\over
\sigma}(1-\rho)}\Subset {\mathcal
U}_1).
\end{equation}  
Derivatives of ${}_\rho\Psi_\sigma$ satisfy, with
universal constant $C$: 
\begin{equation}
\label{Psirhosigma2}
|D^k \left({}_\rho\Psi_\sigma\right)| \leq C^k\left({\sigma \over
1-\rho}\right)^k,
\;k\leq 4.
\end{equation}
uniformly in $\rho, \sigma.$ Of course any
other (fixed) bound for $k$ would do.  

 While it is true that we could just write
$\|\varphi w\|_s\leq c\|\varphi\|_s\|w\|_s,$ for
$s\geq 2,$ to do so would incur at least two derivatives on
$\varphi$ with no gain on $w.$ To avoid this
difficulty, we use the following finer estimates of the $H^2$ norm of
product of functions.

\begin{proposition}
   If $\varphi,\tilde\varphi$ are two smooth, compactly
supported functions with $\tilde\varphi\equiv 1$
on supp $\varphi$ then for every $p\geq 2$
\begin{equation} \label{10} \ \ \ \|\varphi D^pu\|_2\leq C^2\sup_{q\leq
2}\|D^q\varphi\|_{L^{\infty}}\|\tilde\varphi
D^{p-q}u\|_2\ ,\quad {\hbox{and}}\end{equation}
\begin{equation} \label{11} \ \ \ \|\varphi D^pu\|_2\leq C^2\sup_{q\leq
2}\|D^q\varphi\|_{L^{\infty}}\|D^{p-q}u\|_{H^2(supp\
\varphi)}\ .\end{equation}
\end{proposition}

\section{Expanding the norm of the product in (\ref{[Tr,P]})}

The norm of the product in (\ref{[Tr,P]}) will be replaced, as announced,
by the product of the $H^2$ norms, most of which will have as new
functions ${}_{1/m}\Psi_{m''_j}:$ multiplying through by $m,$ 
\begin{multline}\label{onefullpass}{|||{}_0\Psi_m D_t^mu|||_2\over
(m-1)!} 
\lesssim 
\sum_{k=1}^2{\|{}_0\Psi_m ^{(k)}X^2D_t^{m-k}u\|_2\over (m-1)!} +
\ldots 
\\
 +  \sup_{{{m\geq m'-m''\geq 1}
\atop \sum_{j=1}^{m'-m''} m''_j =
m''}\atop {(\sum_{j=1}^{m'-m''} (m''_j+1) =
m'})} 
\left(\prod_{j=1}^{m'-m''}{C\|{}_{1/m}\Psi_{m''_j}D_t^{m''_j}{u'}\|_2\over
m''_j!}\right)
{\|{}_0\Psi_m  \;
X^2D_t^{m-m'}u\|_2\over (m-m'-1)!}
\end{multline}
where, using (\ref{no.ofterms}), the factor $m/m-m'$ which entered on
the right from multiplying through by $m$ and decreasing the last
denominator by one is absorbed in a slightly larger constant
$C^{m'-m''}$ in the product. We have also bounded the terms
$\|h^{(r)}(x,t,u)/r!\|$ by $C^r$ and distributed these constants, one per
term in the product of norms of derivatives of $u'.$  

To unify these two types of terms we could combine them into one sum,
over $k+m'\geq 1,$ but there is nothing new introduced by considering the
couple of extra derivatives which the localizing functions may receive -
there is compensation with decrease in $m$ and we have already seen this
effect - it is essentially one familiar in elliptic regularity proofs by $L^2$
methods, so we will omit the terms with $k>0.$

Now we have seen that we may bring the last localizing function,
${}_0\Psi_m^{(k)},$ out of the last norm and introduce the next function,
${}_{1/m}\Psi_{m-m'},$ identically equal to one on the support of
${}_0\Psi_m,$ with a larger constant
$C_h.$ According to the above Proposition, when bringing a localizing
function out of the norm its $L^\infty$ norm will contribute up to two or,
if already differentiated, perhaps four factors of $m$ with corresponding
decrease in the number of derivatives on $u.$ This disturbs the balance
between number of derivatives and the factorial, but (\ref{no.ofterms}
) shows that
even factors of roughly $(m/m-m)^4$ merely serve to modify the
constant $C_h;$ we conclude that we may pass from one localizing
function to the next without problems. 

That is, applying (\ref{onefullpass}) to its own last term,
with $m$ replaced by $m-m',$ and ignoring $k>0$ for simplicity, we have,
denoting by
${1\over m_2} = 1-(1-{1\over m})(1-{1\over m-m'})$ the band used up by
the supports of the first two localizing functions, which will depend on the
choice of $m',$ and once again ignoring the first term on the right,
\begin{multline}\label{anotherfullpass}{|||{}_{1/m}\Psi_{m-m'}
D_t^{m-m'}u|||_2\over (m-m'-1)!} 
\lesssim 
\\
 +  \sup 
\left(\prod_{j=1}^{\rho'-\rho''}
{C\|{}_{1/ m_2}\Psi_{\rho''_j}D_t^{\rho''_j}{u'}\|_2\over
\rho''_j!}\right)
{\|{}_{1/m}\Psi_{m-m'}  \;
X^2D_t^{m-m'-\rho'}u\|_2\over (m-m'-\rho'-1)!}
\end{multline}
or together, 
\begin{multline}\label{unified2passes}{|||{}_0\Psi_m D_t^mu|||_2\over
(m-1)!} 
\lesssim \sup \left(\prod_{j=1}^{m'-m''}{C\|{}_{1/m}\Psi_{m''_j}D_t^{m''_j}{u'}\|_2\over
m''_j!}\right) \times
\\
\times
\left(\prod_{j=1}^{\rho'-\rho''}
{C\|{}_{1/ m_2}\Psi_{\rho''_j}D_t^{\rho''_j}{u'}\|_2\over
\rho''_j!}\right)
{\|{}_{1/m}\Psi_{m-m'}  \;
X^2D_t^{m-m'-\rho'}u\|_2\over (m-m'-\rho'-1)!}
\end{multline}
where the supremum is over both sets of indices: $m\geq m'+\rho'-m''-\rho''$ and
$\sum_{j=1}^{m'-m''}\sum_{k=1}^{\rho'-\rho''}
(m_j^{''}+\rho_k^{''})= m''+\rho''$
so if we set $s'=m'+\rho'$ and $s''=m''+\rho'',$ we have
a sum over all $m\geq s'-s''\geq 2$ and $\sum_{j+k=2}^{s'-s''}
s_{j+k}''=s''$  while after the first iteration the sum was over all indices such
that $m\geq m'-m''\geq 1,  \sum_{j=1}^{m'-m''} m''_j = m''.$ In both cases,
and for all succeeding ones, the number of such possibilities was seen by
(\ref{partitions}) to be bounded by $C^m.$

We continue this process, pulling the localizing function
${}_{1/m}\Psi_{m-m'}$ out of the last norm and replacing it with
${}_{1/m_2}\Psi_{m-m'-\rho'}$, subjecting that term to the {\it a priori}
estimate, etc. Each time there is a whole `spray' of far lower order terms,
but the number of these is $s'-s'',$ each has a suitable localizing function
which will let us pass to a subsequent one by placing one (universal)
constant with each new copy of $u',$ and in the end we have a product of
on the order of $m$ terms of the form $\|D^ru\|_{H^2(\mathcal{U}_1)}$
all of order $r\leq 4,$ say. (After all, localizing functions need not be
introduced at the last stages - or even in any of the above, until we need to
estimate a given term carefully - for instance, in the product in
(\ref{onefullpass}) the terms could easily have been left as
${\|D_t^{m''_j}u'\|_{H^2(\mathcal{U}_{1/m})}\over  m''_j!}$,) at least
until the time came to subject that term to the {\it a priori} estimate to
reduce its order (in case all other terms had been reduced to lower order). 

We also need to remark at the end that what was true for the first localizing
function, namely (\ref{no.ofterms}), will be a little different, since the next
localizing function may bring not a factor of $m-m'$ with each derivative it
receives but rather the factor (cf. (\ref{Psirhosigma2})) 
$${m-m'\over 1-{1\over m}}= (m-m')\left({m\over m-1}\right)$$
so that, passing from $m-m'$ to $m-m'-n'$ we encounter instead of just 
$${m \over {m-m'}}\leq m'-m''+1$$
an extra factor of $m/m-1,$ possibly to the fourth power; and this may 
keep occurring as the order of the leading term keeps decreasing. For
instance, after a few iterations, the analogous `extra' factors from
(\ref{Psirhosigma2})  will be
$$\left({m\over m-1}\right)\left({m-m_1\over
m-m_1-1}\right)\left({m-m_1-m_2\over m-m_1-m_2-1}\right)\ldots $$ 
or even the fourth power of such a product. But there cannot be more than
$m$ terms in the product and each factor is far less than $2,$ leading to an
easily acceptable constant $C^m$ in the end.

This will prove the bounds for the laft hand side of  (\ref{onefullpass})
$${|||{}_0\Psi_m D_t^mu|||_2\over (m-1)!} \leq C^{m+1}$$
uniformly in $m$ and hence the analyticity of $u$ in $\mathcal{U}_0.$
\qed

\end{document}